\documentstyle[11pt]{article}

\font\bbb=msbm10 scaled\magstep1

\newcommand{\RR}{\mbox{\bbb R}}
\newcommand{\ZZ}{\mbox{\bbb Z}}

\newtheorem{eg}{Example}
\newtheorem{lemma}{Lemma}
\newtheorem{theo}{Theorem}

\newtheorem{conj}{Conjecture}

\begin{document}

\title{\bf A note on the existence of $\{k, k\}$-equivelar polyhedral maps}
\author{{\bf Basudeb Datta} \\
{\em Department of Mathematics, Indian Institute of Science,} \\
{\em Bangalore 560\,012,  India.}  \\
{\small \em e-mail\,: dattab@math.iisc.ernet.in} }
\date{To appear in {\bf `Beitr\"{a}ge zur Algebra und Geometrie}'}

\maketitle

\vspace{-1cm}
\begin{center}
\begin{tabular}{p{5.4in}}
{\small {\bf Abstract.} A polyhedral map is called $\{p,
q\}$-equivelar if each face has $p$ edges and each vertex belongs
to $q$ faces. In \cite{msw2}, it was shown that there exist
infinitely many geometrically realizable $\{p, q\}$-equivelar
polyhedral maps if $q > p = 4$, $p > q = 4$ or $q-3>p =3$. It was
shown in \cite{dn1} that there exist infinitely many $\{4, 4\}$-
and $\{3, 6\}$-equivelar polyhedral maps. In \cite{b}, it was
shown that $\{5, 5\}$- and $\{6, 6\}$-equivelar polyhedral maps
exist. In this note, examples are constructed, to show that
infinitely many self dual $\{k, k\}$-equivelar polyhedral maps
exist for each $k \geq 5$.  Also vertex-minimal non-singular $\{p,
p\}$-pattern are constructed for all odd primes $p$.} \\[2mm]
\end{tabular}
\end{center}

\section{Introduction and results}

A {\em polyhedral complex} (of dimension 2) is collection of
cycles (finite connected 2-regular graphs) together with the
edges and the vertices in the cycles such that the intersection
of any two cycles is empty, a vertex or an edge. The cycles are
called the {\em faces} of the polyhedral complex. For a
polyhedral complex $K$, $V(K)$ denotes its vertex-set and ${\rm
EG}(K)$ denotes its edge-graph or 1-skeleton.  We say $K$ {\em
finite} if $V(K)$ is finite. If ${\rm EG}(K)$ is connected then
$K$ is said to be {\em connected}.

A polyhedral complex is called a {\em polyhedral $2$-manifold}
(or an {\em abstract polyhedron}) if for each vertex $v$ the
faces containing $v$ are of the form $F_1, \dots, F_m$, where
$F_1\cap F_2, \dots, F_{m-1}\cap F_m, F_{m}\cap F_1$ are edges
for some $m\geq 3$. A connected finite polyhedral 2-manifold is
called a {\em polyhedral map}. A {\em combinatorial $2$-manifold}
is a polyhedral 2-manifold whose faces are 3-cycles. A polyhedral
map is called {\em $\{p, q\}$-equivelar} if each face is a
$p$-cycle and each vertex is in $q$ faces. A polyhedral map is
called {\em equivelar} if it is $\{p, q\}$-equivelar for some
$p$, $q$ (cf. \cite{ms, bs, bw, msw1}).

To each polyhedral complex $K$, we associate a pure 2-dimensional
simplicial complex $B(K)$ (called the {\em barycentric
subdivision} of $K$) whose 2-faces are of the form $ueF$, where
$(u, e, F)$ is a flag (i.e., $e$ is an edge of the face $F$ and
$u$ is a vertex of $e$) in $K$. The geometric carrier of $B(K)$ is
called the {\em geometric carrier} of $K$ and is denoted by
$|K|$. Clearly, $K$ is a polyhedral $2$-manifold if and only if
$B(K)$ is a combinatorial 2-manifold (equivalently, $|K|$ is a
2-manifold). A polyhedral 2-manifold $K$ is called {\em
orientable} if $|K|$ is orientable.

An {\em isomorphism} between two polyhedral complexes $K$ and $L$
is a bijection $\varphi \colon V(K) \to V(L)$ such that $(v_1,
\dots, v_m)$ is a face of $K$ if and only if $(\varphi(v_1),
\dots, \varphi(v_m))$ is a face of $L$. Two complexes are called
{\em isomorphic} if there is an isomorphism between them. We
identify two isomorphic polyhedral complexes. An isomorphism from
$K$ to it self is called an {\em automorphism} of $K$. The set
$\Gamma(K)$ of automorphisms of $K$ form a group. A polyhedral
2-manifold $K$ is called {\em combinatorially regular} if
$\Gamma(K)$ is transitive on flags (cf. \cite{ms}).

For a polyhedral 2-manifold $K$, consider the polyhedral complex
$\widetilde{K}$ whose vertices are the faces of $K$ and $(F_1,
\dots, F_m)$ is a face of $\widetilde{K}$ if $F_1, \dots, F_m$
have a common vertex and $F_1 \cap F_2, \dots, F_{m-1} \cap F_m,
F_{m} \cap F_1$ are edges. Then $\widetilde{K}$ is a polyhedral
2-manifold and called the {\em dual} of $K$. If $\widetilde{K}$ is
isomorphic to $K$ then $K$ is called {\em self dual}.

A {\em pattern} is an ordered pair $(M,  G)$, where $M$ is a
connected closed surface in some Euclidean space and $ G$ is a
finite graph on $M$ such that each component of $M \setminus
 G$ is simply connected. The closure of each component of $M
\setminus  G$ is called a {\em face} of $(M,  G)$. For a face
$F$, the closed path (in $ G$) consisting of all the edges and
the vertices in $F$ is called the {\em boundary} of $F$. A
pattern $(M,  G)$ is said to be {\em non-singular} if the
boundary of each face is a cycle. A non-singular pattern is said
to be a {\em polyhedral pattern} if the intersection of any two
face is empty, a vertex or an edge. A pattern $(M,  G)$ is called
a {\em $\{p, q\}$-pattern} if each face contains $p$ edges and
the degree of each vertex in $ G$ is $q$ (cf. \cite{eek}).

If $(M,  G)$ is a polyhedral pattern then clearly the boundaries
of the faces of $(M,  G)$ form a polyhedral map. Conversely, for
a polyhedral map $K$, let $M = |K|$ and $ G = {\rm EG}(K)$. Then
$(M,  G)$ is a polyhedral pattern and the faces of $K$ are the
boundaries of the faces of $(M,  G)$. This pattern $(M,  G)$ is
called a {\em geometric realization} of $K$. A geometric
realization $(M,  G)$ (in some $\RR^n$) is called {\em linear} if
each face of $M$ is a convex polygon and no two adjacent faces
(i.e., faces which share a common edge) lie in the same plane. If
a polyhedral map has a linear geometric realization in $\RR^3$
then it is called {\em geometrically realizable}.

If $f_0(K)$, $f_1(K)$ and $f_2(K)$ are the number of vertices,
edges and faces respectively of a polyhedral complex $K$ then the
number $\chi(K) :=f_0(K)- f_1(K) + f_2(K)$ is called the {\it
Euler characteristic} of $K$. Observe that $\chi(B(K)) =
\chi(K)$. If $u$ and $v$ are vertices of a face $F$ and $uv$ is
not an edge of $F$ then $uv$ is called a {\it diagonal}. Clearly,
if $d(K)$ is the number of diagonals of a polyhedral complex $K$
then $d(K) + f_1(K)\leq {f_0(K) \choose 2}$ and in the case of
equality each pair of vertices belong to a face. A polyhedral map
$K$ is called a {\em weakly neighbourly polyhedral map} (in
short, {\em wnp map}) if each pair of vertices belong to a common
face.

We know (cf. \cite{dn1}) that there exists a unique $\{p,
q\}$-equivelar polyhedral map if $(p, q)= (3, 3)$, $(3, 4)$ or
$(4, 3)$ and there are exactly two $\{p, q\}$-equivelar polyhedral
maps  if $(p, q)= (3, 5)$ or $(5, 3)$. In \cite{msw2}, McMullen et
al. constructed infinitely many geometrically realizable $\{p,
q\}$-equivelar polyhedral maps for each $(p, q) \in \{(r, 4)  : r
\geq 5\} \cup \{(4, s)  :  s \geq 5\} \cup \{(3, k)  :  k \geq
7\}$. In \cite{dn1}, it was shown that there exist infinitely
many $\{4, 4\}$- and $\{3, 6\}$-equivelar polyhedral maps. It was
also shown that there are exactly two neighbourly $\{3,
8\}$-equivelar polyhedral maps and there are exactly 14
neighbourly $\{3, 9\}$-equivelar polyhedral maps.

In \cite{c}, Coxeter constructed a geometrically realizable
combinatorially regular infinite polyhedral 2-manifold whose
faces are hexagons and each vertex is in six faces (namely, $\{6,
6\,|\,3\}$). In \cite{gr}, Gr\"{u}nbaum constructed another
combinatorially regular infinite polyhedral 2-manifold of type
$\{6, 6\}$ (namely, $\{6, 6\}_4$) (cf. \cite{ms}). In \cite{g},
Gott constructed a geometrically realizable infinite polyhedral
2-manifold whose faces are pentagons and each vertex is in five
faces. If $K$ is a $\{p, q\}$-equivelar polyhedral map on $n$
vertices then $d(K) = nq(p - 3)/2$ and $f_1(K) = nq/2$.
Therefore, if $K$ is an $n$-vertex $\{p, p\}$-equivelar
polyhedral map then $np(p-3)/2 + np/2 \leq n(n - 1)/2$ and hence
$n \geq (p - 1)^2$. Clearly, equality holds if and only if $K$ is
a wnp map. Let $\alpha(p)$ denote the smallest $n$ such that
there exist an $n$-vertex $\{p, p\}$-equivelar polyhedral map.
Clearly, the 4-vertex 2-sphere (the boundary of a 3-simplex) is
the unique $\{3, 3\}$-equivelar wnp map. In \cite{b}, Brehm
proved that there exist exactly three $\{4, 4\}$-equivelar wnp
maps and constructed the 16-vertex $\{5, 5\}$-equivelar
polyhedral map $M_{5, 16}$ (of Example \ref{kkmaps}). It was
shown in \cite{bdn} that $M_{5, 16}$ is the unique $\{5,
5\}$-equivelar polyhedral map on $16$ vertices. So, $\alpha(k)=
(k - 1)^2$ for $k \leq 5$. In \cite{b}, Brehm also constructed the
26-vertex $\{6, 6\}$-equivelar polyhedral map $M_{6, 26}$ (of
Example \ref{kkmaps}). Here we show\,:

\begin{theo}\hspace{-1.8mm}{\bf .}\label{t1}
For each $m \geq 3$ and $n\geq 0$, there exist a $2(3^{m - 1} +
2n - 1)$-vertex self dual $\{2m - 1, 2m - 1\}$-equivelar
polyhedral map and a $(3^m  + 2n - 1)$-vertex self dual $\{2m,
2m\}$-equivelar polyhedral map.
\end{theo}

Thus $(2m - 2)^2 \leq \alpha(2m - 1) \leq 2(3^{m - 1} - 1)$ and
$(2m - 1)^2 \leq \alpha(2m) \leq 3^m - 1$ for all $m \geq 3$. In
\cite{n}, using a computer, Nilakantan has shown that there does
not exist any 25-vertex $\{6, 6\}$-equivelar polyhedral map. So,
$\alpha(6) = 26$ and hence there does not exist any $\{6,
6\}$-equivelar wnp map. We believe the following is true\,:

\begin{conj}\hspace{-1.8mm}{\bf .} 
There does not exist any $\{k, k\}$-equivelar wnp map for $k\geq
7$.
\end{conj}

For the existence of an $n$-vertex $\{k, k\}$-pattern $n$ must be
at least $k + 1$. Here we show\,:

\begin{theo}\hspace{-1.8mm}{\bf .} 
There exists a $(p + 1)$-vertex non-singular $\{p, p\}$-pattern
for each odd prime~$p$.
\end{theo}

\section{Examples and proofs of the results}

We first construct infinitely many $\{k, k\}$-equivelar polyhedral
maps. We need these to prove our results. We identify a polyhedral
complex with the set of faces in it.

\begin{eg}\hspace{-1.8mm}{\rm {\bf .}} \label{kkmaps}
{\rm  For $m \geq 3$ and $n \geq 0$, let \vspace{-2mm}
\begin{eqnarray*}
M_{2m -1, 2(3^{m - 1} + 2n - 1)} & \!\!=\!\! & \{F_{i, 2m - 1} : 1
\leq i \leq 2(3^{m - 1} + 2n - 1)\},  \\
M_{2m, 3^m  + 2n - 1} & \!\!=\!\! & \{F_{i, 2m} : 1 \leq i \leq
3^m  + 2n - 1 \},
\end{eqnarray*}
where $b_{2l - 1} = 3^{l - 1} - 1$, $b_{2l} = 2 \times 3^{l-1} -
1$, for $l \geq 1$ and
\begin{eqnarray*}
F_{i, 2m-1} &\!\!=\!\!& (i+b_1, i+b_2, \dots, i+b_{2m-3},
i+b_{2m-2} +n, i+b_{2m-1}+2n), \\
F_{i, 2m} &\!\!=\!\!& (i+b_1, i+b_2, \dots, i+b_{2m-2},
i+b_{2m-1}, i + b_{2m} + n)
\end{eqnarray*}
are cycles ($(2m-1)$-cycles and $(2m)$-cycles respectively) with
vertices from $\ZZ_{2(3^{m - 1} + 2n - 1)}$ and $\ZZ_{3^m  + 2n -
1}$ respectively. Clearly, there are $2m - 1$ faces through each
vertex in  \linebreak $M_{2m - 1, 2(3^{m - 1} + 2n - 1)}$ and
there are $2m$ faces through each vertex in $M_{2m, 3^m  + 2n -
1}$. So, $f_1(M_{2m - 1, 2(3^{m - 1} + 2n - 1)}) =(3^{m-1} + 2n -
1)(2m - 1)$ and $f_1(M_{2m, 3^m  + 2n - 1}) =(3^m + 2n - 1)m$.
Thus, $\chi(M_{2m - 1, 2(3^{m - 1} + 2n - 1)})$ $ = (3^{m-1} + 2n
- 1)(5 - 2m)$ and $\chi(M_{2m, 3^m  + 2n - 1}) = (3^m + 2n - 1)(2
- m)$. By Lemma \ref{l2} below, $M_{2m + 1, 2(3^{m - 1} + 2n -
1)}$ and $M_{2m, 3^m  + 2n - 1}$ are polyhedral maps.  But, by
Lemma \ref{l4}, none of these polyhedral maps are combinatorially
regular.}
\end{eg}

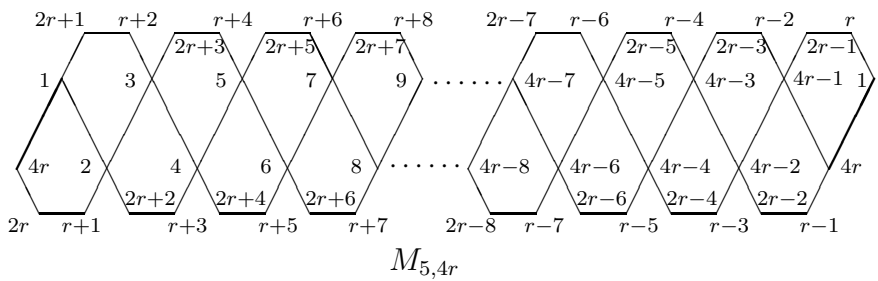
\setlength{\unitlength}{3mm}
\begin{picture}(45,12)(0,-3)

\thinlines

\put(2,0){\line(1,0){2}} \put(6,0){\line(1,0){2}}
\put(10,0){\line(1,0){2}} \put(14,0){\line(1,0){2}}
\put(22,0){\line(1,0){2}} \put(26,0){\line(1,0){2}}
\put(30,0){\line(1,0){2}} \put(34,0){\line(1,0){2}}

\put(4,8){\line(1,0){2}} \put(8,8){\line(1,0){2}}
\put(12,8){\line(1,0){2}} \put(16,8){\line(1,0){2}}
\put(24,8){\line(1,0){2}} \put(28,8){\line(1,0){2}}
\put(32,8){\line(1,0){2}} \put(36,8){\line(1,0){2}}

\put(1,2){\line(1,2){3}} \put(4,0){\line(1,2){4}}
\put(8,0){\line(1,2){4}} \put(12,0){\line(1,2){4}}
\put(16,0){\line(1,2){3}} \put(21,2){\line(1,2){3}}
\put(24,0){\line(1,2){4}} \put(28,0){\line(1,2){4}}
\put(32,0){\line(1,2){4}} \put(36,0){\line(1,2){3}}

\put(2,0){\line(-1,2){1}} \put(6,0){\line(-1,2){3}}
\put(10,0){\line(-1,2){4}} \put(14,0){\line(-1,2){4}}
\put(17,2){\line(-1,2){3}} \put(19,6){\line(-1,2){1}}
\put(22,0){\line(-1,2){1}} \put(26,0){\line(-1,2){3}}
\put(30,0){\line(-1,2){4}} \put(34,0){\line(-1,2){4}}
\put(37,2){\line(-1,2){3}} \put(39,6){\line(-1,2){1}}
\put(14,8){\line(1,-2){1}}

\thicklines

\put(1,2){\line(1,2){2}} \put(37,2){\line(1,2){2}}

\put(17.5,1.8){$\cdots\cdots$} \put(19.4,5.5){$\cdots\cdots$}

\put(0.6,-0.6){\mbox{$_{2r}$}} \put(3,-0.6){\mbox{$_{r+1}$}}
\put(5.8,0.6){\mbox{$_{2r+2}$}} \put(7.7,-0.6){\mbox{$_{r+3}$}}
\put(9.8,0.6){\mbox{$_{2r+4}$}} \put(11.7,-0.6){\mbox{$_{r+5}$}}
\put(13.8,0.6){\mbox{$_{2r+6}$}} \put(15.7,-0.6){\mbox{$_{r+7}$}}
\put(20,-0.6){\mbox{$_{2r-8}$}} \put(23.5,-0.6){\mbox{$_{r-7}$}}
\put(25.8,0.6){\mbox{$_{2r-6}$}} \put(27.7,-0.6){\mbox{$_{r-5}$}}
\put(29.8,0.6){\mbox{$_{2r-4}$}} \put(31.7,-0.6){\mbox{$_{r-3}$}}
\put(33.8,0.6){\mbox{$_{2r-2}$}} \put(35.7,-0.6){\mbox{$_{r-1}$}}

\put(1.5,2){\mbox{$_{4r}$}} \put(3.8,2){\mbox{$_2$}}
\put(7.8,2){\mbox{$_4$}} \put(11.8,2){\mbox{$_6$}}
\put(15.8,2){\mbox{$_8$}} \put(21.5,2){\mbox{$_{4r-8}$}}
\put(25.5,2){\mbox{$_{4r-6}$}} \put(29.5,2){\mbox{$_{4r-4}$}}
\put(33.5,2){\mbox{$_{4r-2}$}} \put(37.5,2){\mbox{$_{4r}$}}

\put(2,5.8){\mbox{$_1$}} \put(5.8,5.8){\mbox{$_3$}}
\put(9.8,5.8){\mbox{$_5$}} \put(13.8,5.8){\mbox{$_7$}}
\put(17.8,5.8){\mbox{$_9$}} \put(23.5,5.8){\mbox{$_{4r-7}$}}
\put(27.5,5.8){\mbox{$_{4r-5}$}} \put(31.5,5.8){\mbox{$_{4r-3}$}}
\put(35.4,5.9){\mbox{$_{4r-1}$}} \put(38.2,5.8){\mbox{$_1$}}

\put(1.8,8.5){\mbox{$_{2r+1}$}} \put(5.5,8.5){\mbox{$_{r+2}$}}
\put(8,7.3){\mbox{$_{2r+3}$}} \put(9.7,8.5){\mbox{$_{r+4}$}}
\put(12,7.3){\mbox{$_{2r+5}$}} \put(13.7,8.5){\mbox{$_{r+6}$}}
\put(16,7.3){\mbox{$_{2r+7}$}} \put(17.7,8.5){\mbox{$_{r+8}$}}
\put(21.8,8.5){\mbox{$_{2r-7}$}} \put(25.5,8.5){\mbox{$_{r-6}$}}
\put(28,7.3){\mbox{$_{2r-5}$}} \put(29.7,8.5){\mbox{$_{r-4}$}}
\put(32,7.3){\mbox{$_{2r-3}$}} \put(33.7,8.5){\mbox{$_{r-2}$}}
\put(36,7.3){\mbox{$_{2r-1}$}} \put(37.7,8.5){\mbox{$_{r}$}}

{\large \put(17.5,-2.5){\mbox{$M_{5, 4r} $}} }
\end{picture}

\begin{lemma}\hspace{-1.8mm}{\rm {\bf .}} \label{l1}
For a collection ${\cal C}$ of cycles, let ${\bf \bar{\cal C}}$
be the $2$-dimensional pure simplicial complex whose $2$-faces
are of the form $xyF$, where $F \in {\cal C}$ and $xy$ is an edge
in $F$. If $B({\cal C})$ is as defined earlier then the following
three are equivalent.
\begin{enumerate}
   \item[\mbox{\rm (i)}] $B({\cal C})$ is a combinatorial $2$-manifold.
   \item[\mbox{\rm (ii)}] ${\bf \bar{\cal C}}$ is a combinatorial $2$-manifold.
   \item[\mbox{\rm (iii)}] For any vertex $v$, the cycles
   containing $v$ are of the form $F_1 = (v, v_{1, 1}, \dots, v_{1,
n_1}), \dots$, $F_m = (v, v_{m, 1}, \dots, v_{m, n_m})$ such that
$v_{1, n_1} = v_{2, 1}, \dots, v_{m - 1, n_{m - 1}} = v_{m, 1}$,
$v_{m, n_m} = v_{1, 1}$ for some $m \geq 2$.
   \end{enumerate}
\end{lemma}

\noindent {\bf Proof.} Clearly, $B({\cal C})$ is a subdivision of
${\bf \bar{\cal C}}$. Therefore, (i) and (ii) are equivalent.

For a 2-dimensional pure simplicial complex $X$, the link of a
vertex $v$ is the graph ${\rm lk}_X(v)$ whose vertex-set is $\{u
\in V(X) ~ : ~ uv \in X\}$ and edge-set is $\{xy ~ : ~ xyv \in
X\}$. Clearly, $X$ is a combinatorial 2-manifold if and only if
${\rm lk}_X(v)$ is a cycle for each $v \in V(X)$.

Let $v$ be a vertex of ${\bf \bar{\cal C}}$. If $v=F\in {\cal C}$
then ${\rm lk}_{\bf \bar{\cal C}}(v)$ is $F$ itself. Let $v$ be a
vertex of ${\bf \bar{\cal C}}$ which is not a cycle in ${\cal C}$.
If the cycles containing $v$ are of the form $F_1 = (v, v_{1, 1},
\dots, v_{1, n_1}), \dots, F_m = (v, v_{m, 1}, \dots, v_{m,
n_m})$ such that $v_{1, n_1} = v_{2, 1}, \dots, v_{m - 1, n_{m -
1}} = v_{m, 1}$, $v_{m, n_m} = v_{1, 1}$ for some $m \geq 2$ then
${\rm lk}_{\bf \bar{\cal C}}(v)$ is the cycle $v_{1, 1} F_1 v_{2,
1} F_2 \cdots v_{m, 1} F_m$. Conversely, if ${\rm lk}_{\bf
\bar{\cal C}}(v)$ is a cycle then, from the definition of ${\bf
\bar{\cal C}}$, ${\rm lk}_{\bf \bar{\cal C}}(v)$ must be of the
form $v_{1, 1} F_1 v_{2, 1} F_2 \cdots v_{m, 1} F_m$, where $F_1
= (v, v_{1, 1}, \dots, v_{1, n_1}), \dots, F_m = (v, v_{m, 1},
\dots, v_{m, n_m})$ such that $v_{1, n_1} = v_{2, 1}, \dots, v_{m
- 1, n_{m - 1}} = v_{m, 1}$, $v_{m, n_m} = v_{1, 1}$. This proves
that (ii) and (iii) are equivalent. \hfill $\Box$

\begin{lemma}\hspace{-1.8mm}{\rm {\bf .}} \label{l2}
$M_{2m - 1, 2(3^{m - 1} + 2n - 1)}$ and $M_{2m, 3^m  + 2n - 1}$
are polyhedral maps for $m \geq 3$, $n \geq 0$.
\end{lemma}

\noindent {\bf Proof\,:} Since $\{i, i+1\}$ is an edge in $M_{2m
- 1, 2(3^{m - 1}+ 2n - 1)}$ for each $i$, ${\rm EG}(M_{2m - 1,
2(3^{m - 1} + 2n - 1)})$ is connected. Similarly, ${\rm
EG}(M_{2m, 3^m  + 2n - 1})$ is connected.

Observe that the faces in $M_{2m - 1, 2(3^{m - 1} + 2n - 1)}$
containing $i$ are $F_i, F_{i-b_2}, F_{i-b_3}, F_{i-b_4}, \dots$,
$F_{i-b_{2m -3}}, F_{i-b_{2m-2}-n}, F_{i-b_{2m-1}-2n}$, where
$F_i = F_{i, 2m-1} = (i+b_1, i+b_2, \dots, i+b_{2m-3}, i+b_{2m-2}
+n, i+b_{2m-1}+2n)$. Clearly, $F_i\cap F_{i - b_3} = \cdots =
F_{i} \cap F_{i-b_{2m-2}-n} = \cdots = F_{i-b_{2m-1}-2n}\cap
F_{i-b_2} = \cdots = F_{i-b_{2m-1}-2n} \cap F_{i-b_{2m-3}} =
\{i\}$.

Since $b_{2j+1} = 2b_{2j} - b_{2j-1}$ for all $j$, $ F_{i -
b_{2l-1}} \cap F_{i-b_{2l}} $ is the edge $\{i, i+b_{2l}-
b_{2l-1}\}$, $F_{i-b_{2l}}\cap F_{i - b_{2l+1}} $ is the edge $
\{i+b_{2l}-b_{2l+1}, i\}$ for $1\leq l\leq m-2$, $F_{i-b_{2m-3}}
\cap F_{i-b_{2m-2}-n}$ is the edge $\{i, i+b_{2m-1}-b_{2m- 2}
+n\}$ and $F_{i-b_{2m-2}-n} \cap F_{i- b_{2m-1} - 2n}$ is the
edge $\{i + b_{2m-3}-b_{2m-2}-n, i\}$. Again, since $2b_{2m-1} +
4n \equiv 0$ (mod $2(3^{m - 1} + 2n - 1)$), $F_{i- b_{2m-1} -
2n}\cap F_{i}$ is the edge $ \{i, i + b_{2m-1} + 2n\}$. Thus, any
pair of faces containing $i$ intersect in either at $i$ or on an
edge through $i$ and the faces containing $i$ form a single cycle
of adjacent faces (sharing a common edge). Therefore, $M_{2m - 1,
2(3^{m - 1} + 2n - 1)}$ is a polyhedral map.

The faces in $M_{2m, 3^m  + 2n - 1}$ containing $i$ are $C_i,
C_{i-b_2}, C_{i-b_3}, \dots, C_{i-b_{2m -1}}, C_{i -b_{2m}- n}$,
\linebreak where $C_i = F_{i, 2m} = (i+b_1, i+b_2, \dots,
i+b_{2m-1}, i+b_{2m} + n)$ and $C_i \cap C_{i - b_3} = \cdots =
C_{i} \cap C_{i-b_{2m-1}} = \cdots = C_{i-b_{2m}-n} \cap C_{i-b_2}
= \cdots = C_{i-b_{2m}-n} \cap C_{i-b_{2m-2}} = \{i\}$. Also,
since $2b_{2m} - b_{2m-1} + 2n \equiv 0$ (mod $3^{m} + 2n - 1$),
$C_{i - b_{2l-1}} \cap C_{i-b_{2l}} $ is the edge $\{i, i+b_{2l}-
b_{2l-1}\}$, $C_{i-b_{2l}}\cap C_{i - b_{2l+1}} $ is the edge $
\{i+b_{2l}-b_{2l+1}, i\}$ for $1\leq l\leq m-1$, $C_{i-b_{2m-1}}
\cap C_{i-b_{2m}-n}$ is the edge $\{i, i - b_{2m} - n\}$ and
$C_{i-b_{2m}-n} \cap C_{i}$ is the edge $ \{i + b_{2m} + n, i\}$.
Thus, any pair of faces containing $i$ intersect in either at $i$
or on an edge through $i$ and the faces containing $i$ form a
single cycle of adjacent faces. Therefore, $M_{2m, 3^{m} + 2n -
1}$ is a polyhedral map. \hfill $\Box$

\bigskip

From the uniqueness of 16-vertex $\{5, 5\}$-equivelar polyhedral
map it follows that $M_{5, 16}$ is self dual. Here we prove.

\begin{lemma}\hspace{-1.8mm}{\rm {\bf .}} \label{l3}
$M_{2m - 1, 2(3^{m - 1} + 2n - 1)}$ and $M_{2m, 3^m + 2n - 1}$ are
self dual for $m\geq 3$ and $n\geq 0$.
\end{lemma}

\noindent {\bf Proof.} Let $\varphi \colon M_{2m - 1, 2(3^{m - 1}
+ 2n - 1)} \to \widetilde{M}_{2m - 1, 2(3^{m - 1} + 2n - 1)}$ be
the mapping given by $\varphi(i) = F_i := F_{-i, 2m-1}$. Clearly
$\varphi$ is a bijection. Consider the face $F_{i}= (i+b_1,
\dots, i+b_{2m-3}, i+b_{2m-2} +n, i+b_{2m-1}+2n)$. Now,
$(\varphi(i + b_1), \dots, \varphi(i + b_{2m-3}), \varphi(i +
b_{2m - 2} + n), \varphi(i + b_{2m - 1} + 2n)) = (F_{- i - b_1},
\dots, F_{-i - b_{2m -3}}, F_{-i-b_{2m-2}-n}, F_{-i-b_{2m-1}-2n})
= \widetilde{\em F}_{-i}$ (say). From the proof of Lemma
\ref{l2}, $\widetilde{F}_{-i}$ is a cycle of adjacent faces
(sharing a common edge) containing the common vertex $-i$.
Therefore, by the definition, $\widetilde{F}_{-i}$ is a face of
$\widetilde{M}_{2m - 1, 2(3^{m - 1} + 2n - 1)}$. This implies
that $\widetilde{M}_{2m - 1, 2(3^{m - 1} + 2n - 1)}$ is
isomorphic to  ${M}_{2m - 1, 2(3^{m - 1} + 2n - 1)}$. Similarly,
$\psi \colon  M_{2m, 3^{m} + 2n - 1} \to \widetilde{M}_{2m, 3^{m}
+ 2n - 1}$, given by $\psi(i) = F_{-i,2m}$ defines an
isomorphism. \hfill $\Box$

\bigskip

Clearly, $\Gamma(M_{2m - 1, 2(3^{m - 1} + 2n - 1)})$ and
$\Gamma(M_{2m, 3^m + 2n - 1})$ are transitive on the vertices and
the faces. Here we prove.

\begin{lemma}\hspace{-1.8mm}{\rm {\bf .}} \label{l4}
$M_{2m - 1, 2(3^{m - 1} + 2n - 1)}$ and $M_{2m, 3^m + 2n - 1}$ are
not combinatorially  regular for all $m\geq 3$ and $n\geq 0$.
\end{lemma}

\noindent {\bf Proof.} Let $\mu = 2(3^{m - 1} + 2n - 1)$. If
$m>3$ then consider the flags ${\cal F}_1 = (0, \{0, b_m - b_{m +
1}\}, F_{- b_{m + 1}})$ and ${\cal F}_2 = (0, \{0, b_{m + 2} -
b_{m + 1}\}, F_{- b_{m + 1}})$ in $M_{2m - 1, \mu}$. If possible
let there exist $\varphi \in \Gamma(M_{2m - 1, \mu})$ such that
$\varphi({\cal F}_1)= {\cal F}_2$. Then $\varphi(0) = 0$,
$\varphi(F_{- b_{m + 1}}) = F_{- b_{m + 1}}$ and hence $\varphi(1
- b_{m + 1}) = - b_{m + 1}$ and $\varphi(1) = 1$. If $m>5$ then,
by considering the faces containing 1, $\varphi(F_{1 - b_{m +
2}}) = F_{1-b_{m + 2}}$, $\varphi(F_{1 - b_{m + 1}}) = F_{1 -
b_{m + 3}}$. These imply $1+b_4- b_{m+3}= \varphi(1 - b_{m + 1})
= - b_{m+1}$ in $\ZZ_{\mu}$, a contradiction. If $m=5$ then
$\varphi(F_{1 - b_{6}}) = F_{1 - b_{8} - n}$ and hence
$1+b_4-b_{8}-n = \varphi(1-b_6) = - b_{6}$ in $\ZZ_{\mu}$. This
is not possible. If $m=4$ then $\varphi(F_{1 - b_{5}}) = F_{1 -
b_{7} - 2n}$ and hence $1+b_4-b_{7} -2n = \varphi(1-b_{-5}) = -
b_{5}$ in $\ZZ_{\mu}$, a contradiction.

For $m = 3$, if $\psi\in \Gamma(M_{5, \mu})$ such that $\psi((0,
\{0, 3+n\}, F_{-b_{4}-n})) = (0, \{0, 13+ 3n\}, F_{-b_{4}-n})$
then $\psi(12+3n) =11+3n$ and $\psi(F_{1-b_4-n}) = F_1$ and hence
$3 = \psi(12+3n) = 11+3n$ in $\ZZ_{\mu}$. This is also not
possible.

Thus, $M_{2m - 1, 2(3^{m - 1} + 2n - 1)}$ always has a pair of
flags ${\cal F}_1$ and ${\cal F}_2$ such that $\varphi({\cal
F}_1) \neq {\cal F}_2$ for all $\varphi\in \Gamma(M_{2m - 1,
2(3^{m - 1} + 2n - 1)})$. So, $M_{2m - 1, 2(3^{m - 1} + 2n - 1)}$
is not combinatorially regular.

Let $\eta= 3^m + 2n - 1$ and $C_i=F_{i, 2m}$. Consider the flags
${\cal C}_1 = (0, \{0, (- 1)^m (b_{m + 2} - b_{m + 1})\},
C_{-b_{m+1}})$ and ${\cal C}_2 = (0, \{0, (- 1)^m (b_{m + 2} -
b_{m + 1})\}, C_{-b_{m+ 2}})$ in $M_{2m, \eta}$. If $\varphi \in
\Gamma(M_{2m, \eta})$ such that $\varphi({\cal C}_1)= {\cal C}_2$.
Then $\varphi({\cal C}_2) = {\cal C}_1$, $\varphi(1 - b_{m + 2})
= - b_{m+1}$ and $\varphi(1) = 1$. If $m>3$ then $\varphi(C_{1 -
b_{m + 2}}) = C_{1 - b_{m + 3}}$ and hence $1+b_2- b_{m+3}=
\varphi(1 - b_{m + 2}) = - b_{m+1}$ in $\ZZ_{\eta}$, a
contradiction. If $m =3$ then $\varphi(C_{1 - b_{5}}) = C_{1 -
b_{6}-n}$ and hence $1+b_2- b_{6}+n= \varphi(1 - b_{5}) = -
b_{4}$ in $\ZZ_{\eta}$. This is not possible. Therefore, by
similar argument as before, $M_{2m, 3^m + 2n - 1}$ is not
combinatorially regular. \hfill $\Box$

\begin{eg}\hspace{-1.8mm}{\rm {\bf .}} {\rm Let ${\cal C}_4$ be
the collection of $4$-cycles of the complete graph $K_5$ on the
vertex set $\ZZ_4\cup\{u\}$ given by ${\cal C}_4 = \{(0, 1, 2,
3), (u, i, i + 1, i + 3) ~ : i\in \ZZ_4\}$. Then $|{\bf \bar{\cal
C}}_4|$ is the torus and hence $(|{\bf \bar{\cal C}}_4|, K_5)$ is
a non-singular $\{4, 4\}$-pattern.}
\end{eg}

\begin{lemma}\hspace{-1.8mm}{\rm {\bf .}} \label{l5}
Suppose ${\cal C}(\pi_p) = \{(0, 1, \dots, p-1),
   (u, i + \pi_p(1), \dots, i + \pi_p({p - 1})) ~ : ~ i \in \ZZ_p\}$
   is a collection of cycles of the complete graph $K_{p + 1}$ on
   the vertex set $\ZZ_p \cup \{u\}$, where $p$ is an odd prime
   and $\pi_p$ is a  permutation of $\ZZ_p\setminus\{0\} =\{1, \dots,
   p-1\}$. If \vspace{-1mm}
   \begin{enumerate}
   \item[{\rm (pp1)}] $\pi_p(i) + \pi_p({p - i}) = p$ for $1 \leq
   i \leq p - 1$,
   \item[{\rm (pp2)}] $\pi_p(\frac{p - 1}{2}) = \frac{p-1}{2}$ and
   \item[{\rm (pp3)}] exactly one of $j$, $- j$ is in $\{\pi_p(2)
   - \pi_p(1), \pi_p(3) - \pi_p(2), \dots, \pi_p(\frac{p + 1}{2})
   - \pi_p(\frac{p - 1}{2})\}$
   \end{enumerate}
   then ${\bf \bar{\cal C}}(\pi_p)$ is a connected combinatorial
   $2$-manifold.
\end{lemma}

\noindent {\bf Proof.} Since edges of cycles of ${\cal C}(\pi_p)$
form a connected graph, ${\rm EG}({\bf \bar{\cal C}}(\pi_p))$ is
connected.

Let $a_i = \pi_p(i + 1)- \pi_p(i)$ for $1 \leq i \leq p - 2$.
Then, by (pp1), $a_i = a_{p - 1 - i}$. Let $r = \frac{p - 3}{2}$.
Then, by (pp1), (pp2), $a_{r + 1} = 1$ and, by (pp3), $\{a_1,
\dots, a_{r + 1}, -a_1, \dots, -a_{r + 1}\} = \ZZ_p \setminus
\{0\}$.

If $r$ is even then the cycles containing $i$ are $(i, u, \dots, i
+ a_1)$, $(i, i + a_1, \dots, i - a_2)$, $(i, i - a_2, \dots, i +
a_3), \dots, (i, i + a_{r - 1}, \dots, i - a_r)$, $(i, i - a_r,
\dots, i + 1)$, $(i, i + 1, i + 2, \dots, i + p - 1)$, $(i, i + p
- 1, \dots, i + a_{r + 2}), \dots, (i, i + a_{2r}, \dots, i-a_{2r
+ 1})$, $(i, i - a_{2r + 1}, \dots, u)$.

If $r$ is odd then the cycles containing $i$ are $(i, u, \dots, i
+ a_1)$, $(i, i + a_1, \dots, i - a_2)$, $(i, i - a_2, \dots, i +
a_3), \dots, (i, i - a_{r - 1}, \dots, i + a_r)$, $(i, i + a_r,
\dots, i + p - 1)$, $(i, i + p - 1, \dots, i + 2, i + 1)$, $(i, i
+ 1, \dots, i - a_{r + 2}), \dots, (i, i + a_{2r}, \dots, i -
a_{2r + 1})$, $(i, i - a_{2r + 1}, \dots, u)$.

The cycles containing $u$ are $(u, \pi_p(1), \dots, \pi_p(p -
1))$, $(u, 1 + \pi_p(1), \dots, 1 + \pi_p(p - 1)), \dots$, $(u, p
- 1 + \pi_p(1), \dots, p - 1 + \pi_p(p - 1))$. Since $\{\pi_p(p -
1), 1 + \pi_p(p - 1), \dots, p - 1 + \pi_p(p - 1)\} = \ZZ_{p}$,
the cycles containing $u$ can be arranged as $(u, \pi_p(i_1),
\dots, \pi_p(j_1)), \dots, (u, \pi_p(i_p), \dots, \pi_p(j_p))$,
where $j_1 = i_2, \dots, j_{p - 1} = i_p, j_p = i_1$. The lemma
now follows by Lemma \ref{l1}. \hfill $\Box$

\medskip

Clearly, $\pi_3$ is the identity permutation and ${\cal
C}(\pi_3)$ is the $4$-vertex $2$-sphere $S^{\,2}_4$. Also,
$\chi({\bf \bar{\cal C}}(\pi_p)) = 2(p + 1) - ({p + 1 \choose 2} +
p(p + 1)) + (p + 1)p = (p + 1)(4 - p)/2$. So, if $p = 4k+1$ for
some $k\geq 1$ then $\chi({\bf \bar{\cal C}}(\pi_p))$ is odd and
hence ${\bf \bar{\cal C}}(\pi_p)$ is non-orientable. Here we
prove.

\begin{lemma}\hspace{-1.8mm}{\rm {\bf .}} \label{l6}
${\bf \bar{\cal C}}(\pi_p)$ is non-orientable for $p>3$.
\end{lemma}

\noindent {\bf Proof.} Let $F= (0, 1,\dots, p-1)$ and $F_i = (u,
i+\pi_p(1), \dots, i+\pi_p(p-1))$ for $1\leq i \leq p-1$. We can
choose a $p$-gonal disc (not necessarily convex) in the plane for
each cycle in ${\bf \bar{\cal C}}(\pi_p)$ so that the disc
corresponding to $F_i$ is attached with that for $F$ along the
common edge $\{i+\pi_p(\frac{p-1}{2}), i+\pi_p(\frac{p+ 1}{2})\}$
for each $i$ and there are no other intersections. This gives us a
$p(p-1)$-gonal disc $D(\pi_p)$. Then there are two edges in
$D(\pi_p)$ corresponding to an edge $jk$ ($j, k\in \ZZ_p$,
$-1\neq j-k\neq 1$) in some cycle $F_i$ and they appear in the
same direction (clockwise or anti-clockwise). Since $|{\bf
\bar{\cal C}}(\pi_p)|$ is homeomorphic to the space obtained by
identifying such pairs of edges (and some more) of $D(\pi_p)$,
$|{\bf \bar{\cal C}}(\pi_p)|$ is non-orientable. \hfill $\Box$

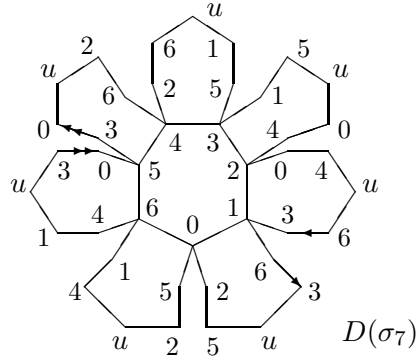
\setlength{\unitlength}{1.8mm}
\begin{picture}(45,26)(-30,2)

\thinlines

\put(8,8){\line(2,3){2}} \put(10,11){\line(2,-1){4}}
\put(10,11){\line(0,1){4}} \put(10,11){\line(-3,-1){3}}
\put(10,15){\line(2,3){2}} \put(10,15){\line(-3,1){3}}
\put(10,15){\line(-3,2){3}} \put(12,18){\line(-3,2){3}}
\put(12,18){\line(-1,3){1}} \put(12,18){\line(1,0){4}}
\put(16,18){\line(1,3){1}} \put(16,18){\line(3,2){3}}
\put(16,18){\line(2,-3){2}} \put(18,15){\line(3,2){3}}
\put(18,15){\line(3,1){3}} \put(18,11){\line(0,1){4}}
\put(18,11){\line(-2,-1){4}} \put(18,11){\line(3,-1){3}}
\put(18,11){\line(2,-3){2}} \put(14,9){\line(1,-3){1}}
\put(14,9){\line(-1,-3){1}}

\put(6,6){\line(1,-1){3}} \put(6,6){\line(1,1){2}}
\put(9,3){\line(1,0){4}} \put(13,3){\line(0,1){3}}
\put(22,6){\line(-1,1){2}} \put(22,6){\line(-1,-1){3}}
\put(15,3){\line(1,0){4}} \put(15,3){\line(0,1){3}}
\put(4,10){\line(1,0){3}} \put(4,16){\line(1,0){3}}
\put(21,10){\line(1,0){3}} \put(21,16){\line(1,0){3}}
\put(4,18){\line(0,1){3}} \put(24,18){\line(0,1){3}}
\put(11,21){\line(0,1){3}} \put(17,21){\line(0,1){3}}
\put(14,26){\line(3,-2){3}} \put(14,26){\line(-3,-2){3}}
\put(2,13){\line(2,3){2}} \put(2,13){\line(2,-3){2}}
\put(26,13){\line(-2,-3){2}} \put(26,13){\line(-2,3){2}}
\put(21,17){\line(3,1){3}} \put(21,17){\line(3,1){3}}
\put(7,17){\line(-3,1){3}}

\put(7,23){\line(-3,-2){3}} \put(7,23){\line(2,-3){2}}
\put(21,23){\line(3,-2){3}} \put(21,23){\line(-2,-3){2}}

\put(20,8){\vector(1,-1){2}} \put(24,10){\vector(-1,0){2}}
\put(4,16){\vector(1,0){2}} \put(4,16){\vector(1,0){2.8}}
\put(7,17){\vector(-3,1){2}} \put(7,17){\vector(-3,1){2.8}}

\put(8,1.4){\mbox{$u$}} \put(12,1){\mbox{\small 2}}
\put(15,1){\mbox{\small 5}} \put(19,1.4){\mbox{$u$}}
\put(0.5,13){\mbox{$u$}} \put(26.5,13){\mbox{$u$}}
\put(2.7,21.5){\mbox{$u$}} \put(24.3,21.5){\mbox{$u$}}
\put(15,26){\mbox{$u$}} \put(4.8,5){\mbox{\small $4$}}
\put(22.5,5){\mbox{\small 3}} \put(11.5,5){\mbox{\small $5$}}
\put(15.7,5){\mbox{\small $2$}} \put(8.5,6.5){\mbox{\small $1$}}
\put(18.5,6.5){\mbox{\small $6$}}  \put(2.5,9){\mbox{\small $1$}}
\put(6.5,10.7){\mbox{\small 4}} \put(13.5,10){\mbox{\small $0$}}
\put(20.5,10.7){\mbox{\small $3$}} \put(24.7,9){\mbox{\small $6$}}
\put(10.5,11.2){\mbox{\small $6$}} \put(16.5,11.2){\mbox{\small
$1$}} \put(4,14){\mbox{\small $3$}} \put(7,14){\mbox{\small $0$}}
\put(10.7,14){\mbox{\small $5$}} \put(16.5,14){\mbox{\small $2$}}
\put(20,14){\mbox{\small $0$}} \put(23,14){\mbox{\small $4$}}
\put(2.5,17){\mbox{\small $0$}} \put(7.5,17){\mbox{\small $3$}}
\put(12.2,16){\mbox{\small $4$}} \put(15,16){\mbox{\small $3$}}
\put(19.4,17){\mbox{\small $4$}} \put(24.7,17){\mbox{\small $0$}}
\put(7.3,19.5){\mbox{\small $6$}} \put(11.8,20){\mbox{\small 2}}
\put(15.2,20){\mbox{\small $5$}} \put(19.8,19.5){\mbox{\small 1}}
\put(5.7,23.2){\mbox{\small $2$}} \put(11.8,23){\mbox{\small $6$}}
\put(15,23){\mbox{\small $1$}} \put(21.7,23){\mbox{\small $5$}}

\put(25,2){\mbox{$D(\sigma_7)$}}
\end{picture}

\vspace{-2mm}

\begin{lemma}$\!\!${\rm {\bf :}} \label{l7} Let $p > 3$  be a prime.
\begin{enumerate}
\item[$(a)$] If $p = 4k + 3$ for some $k \geq 1$ then the
permutation $\sigma_p = (2, 4k + 1)(4, 4k - 1)\cdots (2k, 2k +
3)$ of $\ZZ_p \setminus \{0\}$ satisfies {\rm (pp1)}, {\rm (pp2)}
and {\rm (pp3)} of Lemma $\ref{l5}$.
\item[$(a)$] If $p = 4l + 1$ for some $l \geq 1$ then the
permutation $\rho_p = (1, 4l)(3, 4l - 2)\cdots (2l -1 , 2l + 2)$
of $\ZZ_p \setminus \{0\}$ satisfies {\rm (pp1)}, {\rm (pp2)} and
{\rm (pp3)} of Lemma $\ref{l5}$.
\end{enumerate}
\end{lemma}

\noindent {\bf Proof.} Clearly, $\sigma_p$ and $\rho_p$ satisfy
hypothesis (pp1) and (pp2).

Now, $\{\sigma_p(2) - \sigma_p(1), \dots, \sigma_p(\frac{p+1}{2})
- \sigma_p(\frac{p-1}{2})\} = \{4k, - (4k - 2), 4k - 4,  \dots,
4, - 2, 1\} = \{-2, 4, -6, \dots, -(4k - 2), 4k, -(4k +2)\}$.
Thus $\sigma_p$ satisfies (pp3).

Again, $\{\rho_p(2) - \rho_p(1), \dots, \rho_p(\frac{p+1}{2}) -
\rho_p(\frac{p-1}{2})\}$ $ = \{- (4l - 2), 4l - 4, - (4l - 6),
\dots, 4, - 2, 1\} = \{-2, 4, -6, \dots, (4l - 4), -(4l - 2), -
4l\}$. Thus $\rho_p$ satisfies (pp3). \hfill $\Box$

\medskip

\noindent {\bf Proof of Theorem 1\,:} Let $m \geq 3$ and $n \geq
0$. By Lemma \ref{l2}, $M_{2m - 1, 2(3^{m-1} + 2n - 1)}$ is a
$2(3^{m - 1} + 2n - 1)$-vertex polyhedral map and hence a $\{2m -
1, 2m - 1\}$-equivelar polyhedral map. Again, by Lemma \ref{l2},
$M_{2m, 3^m + 2n - 1}$ is a $(3^m  + 2n - 1)$-vertex polyhedral
map and hence a $\{2m, 2m\}$-equivelar polyhedral map. The
theorem now follows from Lemma \ref{l3}. \hfill $\Box$

\medskip

\noindent {\bf Proof of Theorem 2\,:} Let $p>3$ be a prime and
$K_{p + 1}$ be the complete graph on the vertex set $\ZZ_p \cup
\{u\}$. By Lemma \ref{l7}, there exists a permutation $\pi_p$ of
$\ZZ_p\setminus\{0\}$ which satisfies (pp1), (pp2) and (pp3) of
Lemma \ref{l5}. Let ${\cal C}(\pi_p)$ be as in Lemma \ref{l5}.
Then, by Lemma \ref{l5}, ${\bf \bar{\cal C}}(\pi_p)$ is a
connected combinatorial 2-manifold. So, if $N_p := |{\bf
\bar{\cal C}}(\pi_p)|$ then $(N_p, K_{p + 1})$ is a non-singular
$\{p, p\}$-pattern and the cycles in ${\cal C}(\pi_p)$ are the
boundaries of the faces of $(N_p, K_{p + 1})$. Finally, the
4-vertex 2-sphere $S^{\,2}_4$ gives a $\{3, 3\}$-pattern. This
completes the proof. \hfill $\Box$

\bigskip

\noindent {\bf Acknowledgement\,:} The author is thankful to
Bhaskar Bagchi and Sunanda Bagchi for useful conversations. The
author thanks the anonymous referee for many useful references
and comments which helped to improve the presentation of this
paper.

\end{document}